\documentclass {article}

\usepackage{amsfonts,amsmath,latexsym}
\usepackage{amstext, amssymb, amsthm}
\usepackage[scanall]{psfrag}
\usepackage {graphicx}
\usepackage{verbatim}

\newcommand {\R}{\mathbf{R}}

\newcommand {\B} {\mathbf{B}}

\newcommand {\calC}{\mathcal{C}}

\newtheorem {thm} {Theorem}
\newtheorem* {thm*} {Theorem}

\newtheorem{conj} [thm]{Conjecture}

\begin {document}

\title{ 
{\bf A numerical investigation of level sets of extremal Sobolev functions}}
\author{Stefan Juhnke and Jesse Ratzkin\\
University of Cape Town \\
{\tt juhnke.stefan@gmail.com} and {\tt jesse.ratzkin@uct.ac.za}} 
\maketitle

\begin {abstract} 
In this paper we investigate the level sets of extremal Sobolev functions. For 
$\Omega \subset \R^n$ and $1 \leq p < \frac{2n}{n-2}$, these functions 
extremize the ratio $\frac {\|\nabla u\|_{L^2(\Omega)}}
{\| u \|_{L^p(\Omega)}}$. We conjecture that as $p$ increases the 
extremal functions become more ``peaked" (see the introduction below 
for a more precise statement), and present some numerical evidence to 
support this conjecture. 
\end {abstract}

\section{Introduction}

Let $n \geq 2$, and let $\Omega \subset \R^n$ be a bounded domain 
with piecewise Lipschitz boundary, satisfying a uniform cone 
condition. One can associate a large variety of geometric and 
physical constants to $\Omega$, such as volume, perimeter, 
diameter, inradius, the principal frequency $\lambda(\Omega)$, 
and torsional rigidity $P(\Omega)$ (which is also 
the maximal expected exit time of a standard Brownian particle). 
For more than a century, many mathematicians have investigated how all these quantities 
relate to each other; P\' olya and Szeg\H o's manuscript \cite {PS} provides 
the best introduction to this topic, which remains very active today, 
with many open questions. 

In the present paper we investigate the quantity 
\begin {equation} \label {defn-Cp}
\calC_p(\Omega) = \inf \left \{ \frac{ \int_\Omega |\nabla u|^2 d\mu }
{\left ( \int_\Omega |u|^p d\mu \right )^{2/p} } : u \in W^{1,2}_0(\Omega), 
u \not \equiv 0 \right \}.\end {equation}
The constant $\calC_p(\Omega)$ gives the best constant in the 
Sobolev embedding: 
$$u \in W^{1,2}_0(\Omega) \Rightarrow \| u \|_{L^p(\Omega)} \leq 
\frac{1}{\sqrt{\calC_p(\Omega)}} \| \nabla u \|_{L^2(\Omega)}.$$
By Rellich compactness, the 
infimum in \eqref{defn-Cp} is finite, positive, and realized by an extremal 
function $u^*_p$, which we 
can take to be positive inside $\Omega$ (see, for instance, \cite{GT} 
or \cite {Sau}). The Euler-Lagrange equation for critical points of the 
ratio in \eqref{defn-Cp} is 
\begin {equation} \label {sobolev-pde} 
\Delta u + \Lambda u^{p-1} = 0, \quad \left. u \right |_{\partial \Omega} 
= 0, \end {equation} 
where $\Lambda$ is the Lagrange multiplier. In the case that $u = u_p^*$ is 
an extremal function, a quick integration by parts argument shows that 
the Lagrange multiplier $\Lambda$ is given by 
$$\Lambda = \calC_p(\Omega) \left ( \int_\Omega (u_p^*)^p d\mu \right )^{\frac{2-p}{p}}.$$

It is worth remarking that in two cases the PDE \eqref{sobolev-pde} becomes linear: 
that of $p=1$ and $p=2$. In the case $p=1$, we recover the torsional rigidity as 
$P(\Omega) = (\calC_1(\Omega))^{-1}$, and in the case $p=2$ we recover the 
principal freqency as $\lambda(\Omega) = \calC_2(\Omega)$. These linear problems are 
both very well-studied, from a variety of perspectives, and the literature attached 
to each is huge. From this perspective, the second author and Tom Carroll began a 
research project several years ago, studying the variational problem \eqref{defn-Cp} 
as it interpolates between torsional rigidity and principal frequency, and beyond. (See, for 
instance,  \cite{CR} and \cite{CR2}.) Primarily, we are interested in two central questions: 
\begin {itemize} 
\item Which of the properties of $P(\Omega)$ and $\lambda(\Omega)$ (and their 
extremal functions) also hold for $\calC_p(\Omega)$ (and its extremal functions)? 
\item Can we track the behavior of $\calC_p(\Omega)$ and its extremal function $u_p^*$ 
as $p$ varies? 
\end {itemize} 

Some of our invegistgations have led us conjecture the following. 
\begin {conj} \label{distribution-conj} 
Let $n \geq 2$ and let $\Omega \subset \R^n$ be a bounded domain with piecewise 
Lipschitz boundary satisfying a uniform cone condition. Normalize the 
corresponding (positive) extremal function $u_p^*$ so that
$$\sup_{x \in \Omega} (u_p^*(x)) = 1,$$ 
and define the associated distribution function 
$$\mu_p(t) = | \{ x \in \Omega : u_p^*(x) > t \} |.$$
Then within the allowable range of exponents we have the inequality 
\begin {equation} \label {distribution-ineq} 
1 \leq p < q \Rightarrow \mu_p(t) > \mu_q(t) \quad \textrm{ for almost every }
t \in (0,1). 
\end {equation} 
If $n =2$ the allowable range of exponents is $1 \leq p < q$,and if 
$n \geq 3$ the allowable range of exponents is $1 \leq p < q < \frac{2n}{n-2}$. 
\end {conj} 

Below we will present some compelling numerical evidence in support of this 
conjecture. The remainder of the paper is structured as follows. In 
Section \ref{related-sec} we provide some context for our present 
investigation, and describe some of the related work present in 
the literature. In Section \ref{algorithm-sec} we describe the numerical 
method we use, as well as its theoretical background, and we 
present our numerical results in Section \ref{results-sec}. We 
conclude with a brief discussion of future work and unresolved 
questions in Section \ref{conclusion-sec}. 

{\small {\sc Acknowledgements:} Most of the work described below 
comes from the first author's honors dissertation, completed under the 
direction of the second author. J.~R. is partially supported by the 
National Research Foundation of South Africa. }

\section{Related results} \label{related-sec}

In this section we will highlight some related theorems about principal frequency, 
torsional rigidity, qualitative properties of extremal functions, and other 
quantities. The following is by no means an exhaustive list. 

The distribution function $\mu_p$ is closely related to a variety of rearrangements
of a generic test function $u$ for \eqref{defn-Cp}. One can rearrange the function 
values of a positive function in a variety of ways, and different rearrangements 
will yield different results. One of the most well-used rearrangements is 
Schwarz symmetrization, where one replaces a positive function $u$ on $\Omega$ 
with a radially-symmetric, decreasing function $u^*$ on $B^*$, a ball with the same volume 
as $\Omega$. The rearrangement is defined to be equimeasurable with $u$: 
$$|\{ u > t \}| = |\{ u^* > t \}| \textrm{ for almost every function value }t.$$

Krahn \cite {K} used Schwarz symmetrization 
to prove an inequality conjectured by Rayleigh in 
the late 1880's:  
\begin {equation} \label {faber-krahn1}
\lambda(\Omega) \geq \left ( \frac{|\Omega|} {\omega_n} \right )^{-2/n} \lambda(\B),\end {equation} 
where $\B$ is the unit ball in $\R^n$, and $\omega_n$ its volume. Moreover, equality can only 
occur in \eqref{faber-krahn1} if $\Omega = \B$ apart from a set of 
measure zero. 
In fact, it is straightforward to adapt Krahn's proof to show
\begin {equation} \label {faber-krahn2} 
|\Omega| = |\B| \Rightarrow \calC_p(\Omega) \geq \calC_p(\B),\end {equation} 
with equality occuring if and only if $\Omega = \B$ apart from a set of measure 
zero (see \cite{CR}).
One can also use similar techniques to prove, for instance, that the square has the 
greatest torsional rigidity among all rhombi of the same area \cite{P}.  

However, there is certainly a limit to the results one can prove using 
only Schwarz (or Steiner) symmetrization, and to go further one must apply 
new techniques. Among these, one can rearrange by weighted volume \cite {PW, R, HH}, 
which works well for wedge-shaped domains. One can rearrange by powers 
of $u$, or (more generally) by some function of the level sets of $u$ \cite {PR1, PR2, Tal, Ch}.  
If one is combining domains using Minkowski addition, then the Minkowski 
sup-convolution is a very useful tool \cite {CCS}. 

All these techniques are successful, to varying degrees, when studying \eqref{defn-Cp} 
for a {\bf fixed} value of $p$. However, we are presently at a loss with regards to 
applying them when allowing $p$ to vary. There are comparitively few results comparing 
the behavior of $\calC_p(\Omega)$ and its extremals $u_p^*$ for different 
values of $p$. 

It is well-known \cite {Tru} that as $p \rightarrow \frac{2n}{n-2}$ the solutions $u_p^*$ 
become arbitrarly peaked, and the distribution function $\mu_p(t)$ approaches $0$ on the 
interval $(\epsilon, 1)$ for any $\epsilon > 0$. This behavior is a reflection of the 
fact that the Sobolev embedding is not compact for the critical 
exponent of $\frac{2n}{n-2}$, and the loss of compactness is due to the 
fact that the functional in \eqref{defn-Cp} is invariant under conformal transformation 
for this exponent. Thus, it is interesting to 
understand the asymptotics as $p \rightarrow \frac{2n}{n-2}$. A partial list 
of such results includes an asymptotic expansion of $\calC_p(\Omega)$ due to 
van den Berg \cite{vdB} and a theorem of Flucher and Wei \cite {FW} (see 
also \cite{BF}) determining the asymptotic location of the maximum of the 
extremal $u_p^*$. Additionally, P. L. Lions \cite {L1, L2} started a program to understand 
the loss of compactness, due to concentration of solutions, for a variety of 
geometric problems in functional analysis and PDE. R. Schoen and Y.-Y. Li (among 
others) have 
exploited this concentration-compactness phenomenon to understand the 
problem of prescribing the scalar curvature of a conformally flat metric. 

We remark that until now we had scant evidence for Conjecture \ref{distribution-conj}. 
Namely, we knew in advance that the extremals become arbitrarily peaked 
as $p$ approaches the critical exponent, and we knew that in the very special case 
$\Omega = \B$ we have $\mu_1 (t) > \mu_2 (t)$.

\section{Our numerical algorithm} \label{algorithm-sec} 

Our numerical method is borrowed from foundational work of 
Choi and McKenna \cite{CM} and Li and Zhou \cite {LZ}, and its 
theoretical underpinning is the famous ``mountain pass" method 
of Ambrosetti and Rabinowitz \cite {AR}. Within our range of allowable 
exponents, Rellich compactness exactly implies that the functional 
\eqref{defn-Cp} satisfies the Palais-Smale condition, and so the mountain 
pass theorem of \cite {AR} implies the existence of a minimax critical point. A later 
refinement of Ni \cite {Ni} implies that in fact a minimax critical point 
lies on the Nehari manifold, defined by 
\begin {equation} \label {nehari-m'fold}
\mathcal{M} = \left \{ u \in W^{1,2}_0(\Omega): u \not \equiv 
0, \int_\Omega |\nabla u |^2 - u^p d\mu = 0 \right \}.
\end {equation} 

To find critical points, we project onto $\mathcal{M}$, using 
the operator 
\begin {equation} \label {nehari-proj} 
P_{\mathcal{M}} (u) = \left ( \frac{\int_\Omega |\nabla u|^2 d\mu}
{\int_\Omega |u|^p d\mu } \right ) ^{\frac{1}{p-2}} u. 
\end {equation} 
Our goal will be to find mountain pass critical points of the 
associated functional 
\begin {equation} \label{sobolev-functional}
\mathcal{I} (u) = \int_\Omega \frac{1}{2} |\nabla u|^2 - 
\frac{1}{p} |u|^p d\mu,\end {equation}
which lie on the Nehari manifold defined in \eqref{nehari-m'fold}. 
Observe that the Frechet derivative of $\mathcal{I}$ is 
\begin {eqnarray*}
\mathcal{I}'(u) (v) & = & \left. \frac{d}{d\epsilon} \right|_{\epsilon = 0} 
\mathcal{I} (u + \epsilon v) \\ 
& = & \int_\Omega \langle \nabla u, \nabla v\rangle - u^{p-1} v d\mu,
\end {eqnarray*} 
so that, after integrating by parts, we can find the direction $v$ of steepest descent 
by solving the equation 
\begin {equation} \label {steepest-descent}
2\lambda \Delta v = -\Delta u - u^{p-1}.\end {equation} 
We are free to choose $\lambda>0$ as a normalization constant, 
and choose it so that $\int_\Omega |\nabla v|^2 d\mu = 1$. (It is 
well-known that by the Poincar\'e inequality this $H^1$-norm is 
equivalent to the $W^{1,2}$-norm.) An expansion of the difference 
quotient (using our normalization of $v$) shows 
$$\frac{\mathcal{I}(u+\epsilon v) - \mathcal{I}(u)}{\epsilon} = 
-2\lambda + \mathcal{O}(\epsilon),$$
so choosing $\lambda>0$ does indeed correspond to the direction 
of steepest descent of $\mathcal{I}$, rather than the direction 
of largest increase.

At this point we remark on the importance of taking $p>2$. In the 
superlinear case $u_0\equiv 0$ is a local minimum and, so 
long as $u \not \equiv 0$ we have $\mathcal {I}(ku) < 0$ for $k>0$ 
sufficiently large. Thus, for any path $\gamma(t)$ joining $u_0$ to 
$ku_{\textrm{guess}}$, the function $h_\gamma(t) = \mathcal{I}(\gamma(t))$
will have a maximum at some value $t_\gamma$. We can imagine 
varying the path $\gamma$ and finding the lowest such maximal 
value, which is exactly our mountain pass critical point.  

We will begin with an initial guess $u_{\textrm{guess}}$
which is positive inside $\Omega$ and $0$ on $\partial \Omega$, 
and let $u_1 = P_{\mathcal{M}}(u_{\textrm{guess}})$. Thereafter 
we apply the following algorithm: 
\begin {enumerate} 
\item Given $u_k$, we compute the direction of steepest descent $v_k$
using \eqref{steepest-descent}. 
\item If $\| v_k \|_{W^{1,2}(\Omega)}$ is suffiently small we stop the 
algorithm, and otherwise we let $u_{k+1} = P_{\mathcal{M}}(u_k+v_k)$
\item If $\mathcal{I}(u_{k+1})< \mathcal{I}(u_k)$ then we repeat the 
entire algorithm starting from the first step. Otherwise we replace $v_k$ 
with $\frac{1}{2} v_k$ and recompute $u_{k+1}$. 
\item Upon the completion of this algorithm, we test our numerical solution 
to verify that it does indeed solve the PDE \eqref{sobolev-pde} weakly. 
\end {enumerate}

Several remarks are in order. The algorithm outlined above is exactly the 
one proposed by Li and Zhou in \cite {LZ}. They proved convergence of 
the algorithm under a wide variety of hypotheses, which include the 
superlinear ($p>2$) case of \eqref{defn-Cp} and \eqref{sobolev-functional}. 
However, they do not claim convergence of the algorithm in the 
sublinear case, and in this case the algorithm fails. On the other hand, we 
are able to verify that in the superlinear case the algorithm coverges to 
a positive (weak) solution of the PDE \eqref{sobolev-pde}, so we are
confident we have reliable data in this case. We present this data in the 
next section. 

In this algorithm we must repeatedly solve the linear PDE \eqref{steepest-descent}, 
which we do in the weak sense, using biquadratic (nine-noded) quadrilateral 
finite elements. In each 
of these steps we replace the corresponding integrals with sums over the 
corresponding elements. We outline this numerical step in the paragraphs below. 

In this computation we take $u$ as known at the mesh points (by an initial guess or 
by the result of a previous iteration). Writing $\overline{v} = 2\lambda v + u$, the 
solution to \eqref{steepest-descent} is given by the solution to
\begin{align}
\Delta \overline{v} = -u^{p-1} \label{descentDirEqu2}
\end{align}
from which we can recover the steepest descent direction $v$. 

To solve for $\overline{v} \in W_0^{1,2}(\Omega)$ we will solve the weak form of \eqref{descentDirEqu2}, {\it i.e.}  
\begin{align}
\int_{\Omega}\nabla w(x) \cdot \nabla \overline{v}(x) dx = \int_{\Omega}w(x)u(x)^{p-1}dx
\end{align}
for any test function $w \in W^{1,2}_0(\Omega)$. 
We will now derive the finite element formulation based on the methods presented by Fish and Belytschko \cite{FB}. We firstly notice that we can split up our integral as a sum of the integrals over the individual element domains $\Omega^e$:
\begin{align*}
\sum_{e=1}^{n_{el}}\bigg{\{}\int_{\Omega^e} \nabla w^e(x)\nabla \overline{v}^{e}(x)dx - \int_{\Omega^e} w^{e}(x)(\overline{v}^e(x))^{p-1}dx\bigg{\}} = 0.
\end{align*}
Now we now write our functions $w$ and $\overline{v}$ in terms of their finite element approximations as:
\begin{align*}
w(x)\approx w^{h}(x) = \boldsymbol{N}(x)\boldsymbol{w}, \qquad \overline{v}(x)\approx \overline{v}^{h}(x) = \boldsymbol{N}(x)\boldsymbol{d},
\end{align*}
where $\boldsymbol{N}$ are quadratic shape functions with value 1 at their corresponding mesh point and value 0 at all other mesh points, while $\boldsymbol{w}$ , $\boldsymbol{d}$ are vectors of nodal function values. The gradients of $w$ and $\overline{v}$ can then be written as
\begin{align*}
\nabla w \approx \boldsymbol{B}(x)\boldsymbol{w}, \qquad \nabla \overline{v} \approx \boldsymbol{B}(x)\boldsymbol{d},
\end{align*}
where $\boldsymbol{B}$ are the gradients of the shape functions. We can rewrite the above expressions for the element level as
\begin{align*}
w^{e}(x)\approx \boldsymbol{N}^{e} (x)\boldsymbol{w}^{e}, \quad \overline{v}^{e}(x)\approx \boldsymbol{N}^{e}(x) \boldsymbol{d}^{e}, \quad
\nabla w^{e} \approx \boldsymbol{B}^{e} (x)\boldsymbol{w}^{e}, \quad \nabla \overline{v}^{e} \approx \boldsymbol{B}^{e}(x)\boldsymbol{d}^{e}.
\end{align*}
Rewriting the integral using these approximations leaves us with
\begin{align*}
\sum_{e=1}^{n_{el}}\bigg\{\int_{\Omega^e}\boldsymbol{w}^{e^{T}}\boldsymbol{B}^{e^{T}}(x)\boldsymbol{B}^{e}(x)\boldsymbol{d}^{e}dx - \int_{\Omega^e}\boldsymbol{w}^{e^T}\boldsymbol{N}^{e^{T}}(x)(\boldsymbol{N}^e(x)\boldsymbol{d}^e)^{p-1}dx\bigg{\}} = 0,
\end{align*}
since ${\boldsymbol{B}^{e}(x) \boldsymbol{w}^{e}}^{T} = \boldsymbol{w}^{e^{T}} \boldsymbol{B}^{e^{T}}(x)$ and ${\boldsymbol{N}^{e}(x) \boldsymbol{w}^{e}}^{T} = \boldsymbol{w}^{e^{T}} \boldsymbol{N}^{e^{T}}(x)$.
We notice that we can take the constants $\boldsymbol{w}^{{e}^{T}}$ and $\boldsymbol{d}^{e}$ outside of the integral to give
\begin{align*}
\sum_{e=1}^{n_{el}}\boldsymbol{w}^{e^{T}}\bigg{\{}\int_{\Omega^e}\boldsymbol{B}^{e^{T}}(x)\boldsymbol{B}^{e}(x)dx\boldsymbol{d}^{e} - \int_{\Omega^e}\boldsymbol{N}^{e^{T}}(x)(\boldsymbol{N}^e(x)\boldsymbol{d}^e)^{p-1}dx\bigg{\}} = 0.
\end{align*}
Letting
\begin{align*}
\boldsymbol{K}^{e} = \int_{\Omega^e}\boldsymbol{B}^{e^{T}}(x)\boldsymbol{B}^{e}(x)dx \quad \textrm{and} \quad  \boldsymbol{f}^{e} = \int_{\Omega^e}\boldsymbol{N}^{e^{T}}(x)(\boldsymbol{N}^e(x)\boldsymbol{d}^e)^{p-1}dx
\end{align*}
and using the gather matrix to write
\begin{align*}
\boldsymbol{w}^{e} = \boldsymbol{L}^{e}\boldsymbol{w}, \quad \boldsymbol{d}^{e} = \boldsymbol{L}^{e}\boldsymbol{d},
\end{align*}
we get
\begin{align*}
\boldsymbol{w}^{T}\bigg{(}\sum_{e=1}^{n_{el}}\boldsymbol{L}^{e^{T}}\boldsymbol{K}^{e}\boldsymbol{L}^{e}\boldsymbol{d} - \sum_{e=1}^{n_{el}}\boldsymbol{L}^{e^{T}}\boldsymbol{f}^{e}\bigg{)} = 0.
\end{align*}
Further letting
\begin{align*}
K = \sum_{e=1}^{n_{el}}\boldsymbol{L}^{e^{T}}\boldsymbol{K}^{e}\boldsymbol{L}^{e}\boldsymbol{d} \quad \textrm{and } \boldsymbol{f} = \sum_{e=1}^{n_{el}}\boldsymbol{L}^{e^{T}}\boldsymbol{f}^{e} = 0,
\end{align*}
we end up with
\begin{align*}\boldsymbol{w}^{T}(\boldsymbol{Kd-f}) = 0 \quad \forall \boldsymbol{w}.
\end{align*}
Since we know that $w \in W_0^{1,2}$ is arbitrary we therefore solve the discrete finite element form 
\begin{align}\boldsymbol{Kd} = \boldsymbol{f} \label{Kdf},
\end{align}
with $\boldsymbol{N}\boldsymbol{d}$ the finite element approximation to $\overline{v}$ from which we can recover the steepest descent direction $v$.

\section{Numerical results} \label{results-sec}

In this section we describe our numerical results. We implemented 
the algorithm described in Section~\ref{algorithm-sec} using in 
{\sc Matlab}, and all the figures displayed below come from this implementation. 

\begin {figure}[h]
\begin {center}
\includegraphics[width=4in]{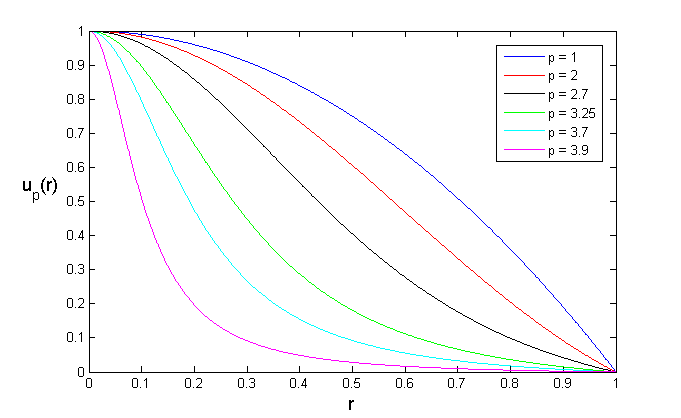}
\caption{Extremal Sobolev functions for a four-dimensional 
unit ball}
\label {4-ball-solns}
\end {center}
\end {figure}

We first implemented our method on a unit ball of dimension four. In this 
case, the solution is radially symmetric, so we only need to solve an ODE. We 
display a plot of these solutions in Figure~\ref{4-ball-solns}.

\begin {figure}[h]
\begin {center}
\includegraphics[width=4in]{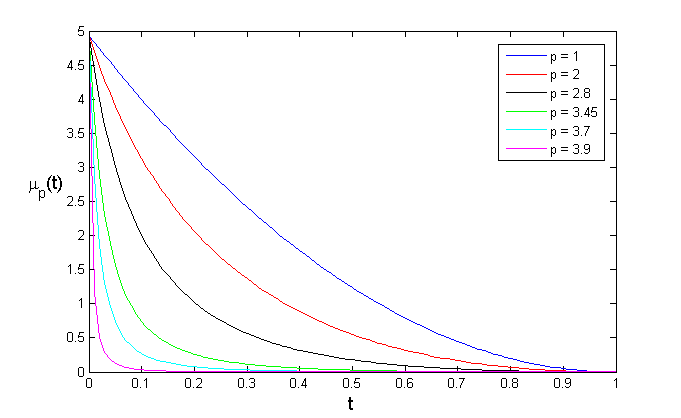}
\caption{Distributions of extremal Sobolev functions for a four-dimensional 
unit ball}
\label {4-ball-distributions}
\end {center}
\end {figure}

We also display a plot of the corresponding distribution functions in 
Figure~\ref{4-ball-distributions}. 

Observe that, as we expected, the distribution function appears to be monotone, and 
that as $p \rightarrow 4 = \frac{2n}{n-2}$ the solution becomes arbitrarily 
concentrated at the origin. 

We can verify that we are indeed finding solutions to the correct PDE. For the case $p=1$
and $p=2$ we can compute the solutions analytically, and verify directly that our numerical 
solution agrees quite well. These are (up to a constant multiple) 
$$u_1^*(r) = 1-r^2, \qquad u_2^*(r) = r^{\frac{2-n}{2}} J_{\frac{n-2}{2}}
\left ( j_{\frac{n-2}{2}} r\right ),$$
where $J_a$ is the Bessel function of the first kind of index $a$ and $j_a$ 
is its first positive zero. For other values of $p$ we can verify that we have found 
a weak solution of \eqref{sobolev-pde}. As the solution is {\it a priori} radial, we 
know that the weak form of the PDE is 
\begin{align}
WT_w(u):= \int_{0}^{1} \bigg[-r^{1-n}\frac{\partial w(r)}{\partial r} \bigg( r^{n-1}\frac{\partial u(r)}{\partial r} \bigg) + w(r)\Lambda u(r)^{p-1}\bigg]r^{n-1}dr & = 0 \label{weakTest}
\end{align}
The above lends itself well to testing via finite element approximation. A random test function $w(r)$ is created by randomly generating numbers at the mesh points and $WT_w(u)$ is evaluated by Gauss quadrature. For comparison purposes, the functions $u$ are normalized so that $\sup(u)$ = 1. This requires that $\Lambda$ be rescaled ($\Lambda$ is set equal to 1 in the algorithm for simplicity), and the appropriate rescaling is then given by $a^{2-p}$ where $a$ is the factor normalizing $u$. This rescaling is derived from the fact that if $u$ solves
\begin{align}
\Delta u + u^{p-1} = 0 \label{thisEqu}
\end{align} 
then $au$ solves 
\begin{align*}
\Delta(au) + a^{2-p}(au)^{p-1} = 0,
\end{align*} 
by simply multiplying \eqref{thisEqu} by $a$.
\\ \\
We generate values of $WT_w(u)$ for a number of test functions $w$ and examine the average magnitude. As alluded to previously, the result of the test \eqref{weakTest} is that for solution candidate functions derived from our algorithm for $2 \leq p < \frac{2n}{n-2}$ and for $p=1$, we have $WT_w(u)$ very close to zero, meaning that we can be confident that we have found appropriate solutions.

\begin {figure}[h]
\begin {center}
\includegraphics[width=4in]{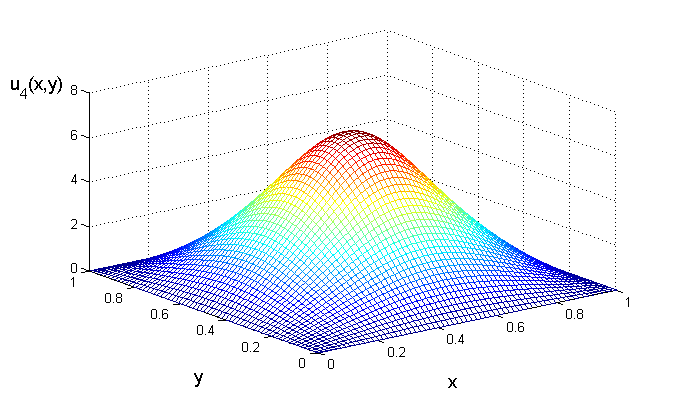}
\caption{Extremal Sobolev function for $p=4$ on a unit square}
\label {square-4}
\end {center}
\end {figure}

\begin {figure}[h]
\begin {center}
\includegraphics[width=4in]{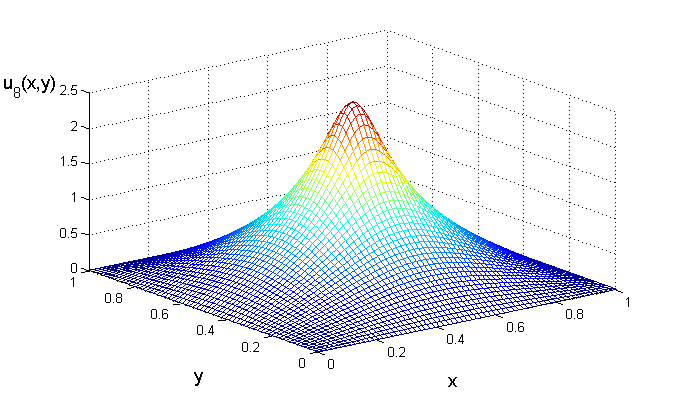}
\caption{Extremal Sobolev function for $p=8$ on a unit square}
\label {square-8}
\end {center}
\end {figure}

\begin {figure}[h]
\begin {center}
\includegraphics[width=4in]{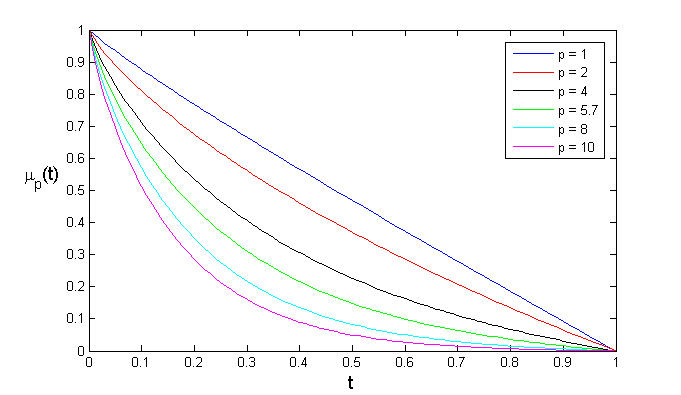}
\caption{Distributions of extremal Sobolev functions for a 
unit square in the plane}
\label {square-distributions}
\end {center}
\end {figure}

Next we implemented our algorithm in a unit square in the plane. We 
plot below our numerical solution for both $p=4$ (Figure~\ref{square-4})
and $p=8$ (Figure~\ref{square-8}), and the distribution function 
for several values of $p$ (Figure~\ref{square-distributions}). 

Again we verify that our numerical algorithm does find a weak solution 
of~\eqref{sobolev-pde}. This time we define
\begin{align}
WT_w(u):= \int_{\Omega}[-\nabla u(x) \nabla w(x) + w\Lambda u(x)^{p-1} ]dx \label{weakTestSquare2}
\end{align}
and again compute $WT_w(u)$ for our candidate solutions, with appropriate rescalings as described previously. We have closely matched the result of Choi and McKenna for the case $p=4$, which means that we should be able to use the value $WT_w(u^*_4)$ as a gauge for how close to zero $WT_w(u)$ should be for appropriate solutions. Again we find that for $2 \leq p < \frac{2n}{2-n}$ and $p=1$ we get values of $WT_w(u)$ very close to zero and of the same magnitude as $WT_w(u^*_4)$.

\begin {figure}[h]
\begin {center}
\includegraphics[width=4in]{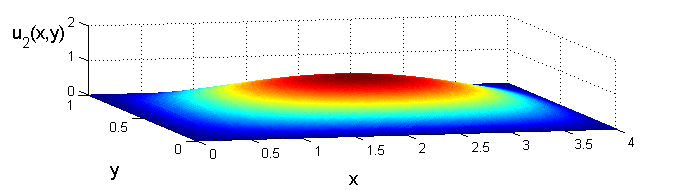}
\caption{Extremal Sobolev function for $p=2$ on a $1
\times 4$ rectangle}
\label {rectangle-2}
\end {center}
\end {figure}

\begin {figure}[h]
\begin {center}
\includegraphics[width=4in]{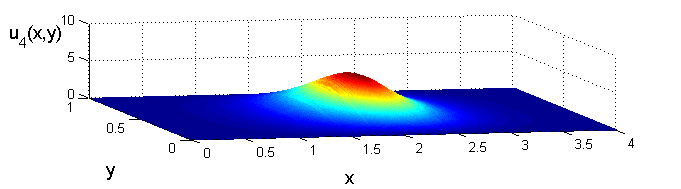}
\caption{Extremal Sobolev function for $p=4$ on a $1\times 
4$ rectangle}
\label {rectangle-4}
\end {center}
\end {figure}

\begin {figure}[h]
\begin {center}
\includegraphics[width=4in]{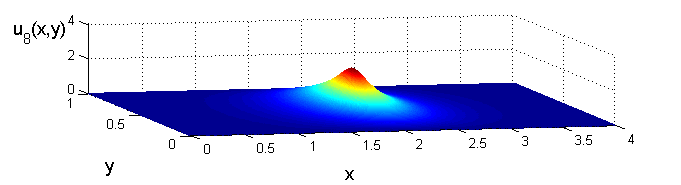}
\caption{Extremal Sobolev function for $p=8$ on a $1
\times 4$ rectangle}
\label {rectangle-8}
\end {center}
\end {figure}

\begin {figure}[h]
\begin {center}
\includegraphics[width=4in]{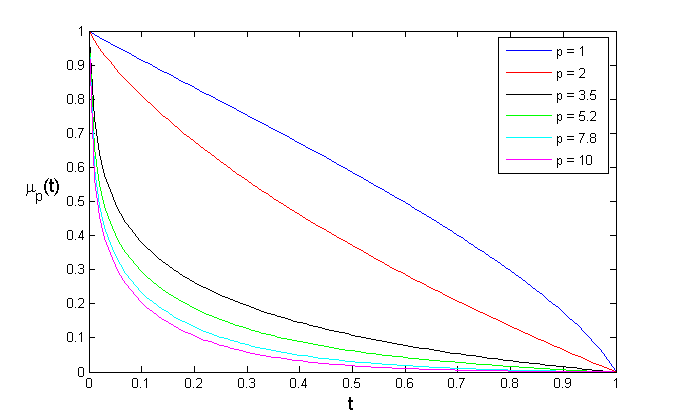}
\caption{Distributions of extremal Sobolev functions for a 
$1\times 4$ rectangle}
\label {rectangle-distributions}
\end {center}
\end {figure}

Finally we implemented our algorithm on a rectangle of width $1$ and length $4$ in the plane. Below 
we plot the our numerical solutions for $p=2$ (Figure~\ref{rectangle-2}), $p=4$ (Figure~\ref{rectangle-4}), 
and $p=8$ (Figure~\ref{rectangle-8}), as well as the distribution function for several values 
of $p$ (Figure~\ref{rectangle-distributions}). We use the same test as we did in the case of the unit square to 
verify that in the case of he $1\times 4$ rectangle we have indeed found (weak) numerical solutions of~\eqref{sobolev-pde}.

\section{Outlook} \label{conclusion-sec} 

The present paper is only the start of our numerical and theoretical investigations
into Conjecture \ref{distribution-conj}. We would like to verify our results on some 
more planar domains, such as triangles and parallelograms. 
Next we anticipate numerical compuations
for higher dimensional objects, such as cubes and parallelpipeds, in the 
super-linear case, as well as possibly some ring domains. We will also need to 
develope a new numerical algorithm which yields reliable results for 
$1 < p < 2$. Finally, we hope that our numerical data provides enough 
insight to rigorously prove our conjecture.

\end{document}